\documentclass[11pt,a4paper,twoside,reqno]{amsart}
\usepackage[T1]{fontenc}  
\usepackage[frenchb]{babel}
\usepackage{amsmath}
\usepackage{array}
\usepackage{amssymb}			
\usepackage{indentfirst}
\usepackage{geometry}
\usepackage{eurosym}
\usepackage{graphicx}
\usepackage{amsthm}
\usepackage{paralist}
\usepackage{amscd}
\usepackage{boxedminipage}
\usepackage{cmbright}
\usepackage{hyperref}
\usepackage{verbatim}
\usepackage{amsbsy}     
\usepackage{setspace}   

\newcommand {\C}{\mathbb{C}}

\newcommand {\N}{\mathbb{N}}

\newcommand {\Q}{\mathbb{Q}}
\newcommand {\R}{\mathbb{R}}	
\newcommand {\Z}{\mathbb{Z}}

\newcommand {\CC}{\mathcal{C}}

\newcommand {\MM}{\mathcal{M}}

\renewcommand {\line}{\vspace{0.25cm}}

\newcommand {\Vect}{\text{\textup{\textrm{Vect}}}}

\renewcommand {\l}{\left}
\renewcommand {\r}{\right}
\newcommand {\ve}{\varepsilon}

\newcommand {\vphi}{\varphi}
\newcommand {\rg}{\text{rg}}

\newtheoremstyle{def}{}{}{}{}{}{ :}{\newline}{}\theoremstyle{def}
\newtheorem{definition}{D\'efinition }[section]

\newtheorem{notation}{Notation }[definition]

\newtheoremstyle{theo}{}{}{}{}{\bfseries}{ :}{\newline}{}\theoremstyle{theo}
\newtheorem{theoreme}{Th\'eor\`eme }[section]

\newtheorem{proposition}[theoreme]{Proposition }

\newtheorem{conjecture}{Conjecture }[section]
\newtheorem{lemme}[theoreme]{Lemme }

\newtheorem{corollaire}[theoreme]{Corollaire }

\theoremstyle{remark}

\newenvironment{theo}
	{\line\line\hspace{-0.5cm}\begin{tabular}{!{\vrule width 2pt}l}\begin{minipage}[t]{14.4cm}\begin{theoreme}}
	{\end{theoreme}\end{minipage}\end{tabular}\line}
\newenvironment{prop}
	{\line\line\hspace{-0.5cm}\begin{tabular}{!{\vrule width 2pt}l}\begin{minipage}[t]{14.4cm}\begin{proposition}}
	{\end{proposition}\end{minipage}\end{tabular}\line}
\newenvironment{lem}
	{\line\line\hspace{-0.5cm}\begin{tabular}{!{\vrule width 2pt}l}\begin{minipage}[t]{14.4cm}\begin{lemme}}
	{\end{lemme}\end{minipage}\end{tabular}\line}

\begin{document}

\title{{\large G\'en\'eralisations quantitatives du crit\`ere d'ind\'ependance lin\'eaire de Nesterenko} \\ {\small A l'occasion de la conf\'erence THUE 150 en l'honneur du 150\`eme anniversaire d'Axel Thue}}

\author{Simon Dauguet}
\address{Universit\'e Paris-Sud, Laboratoire de Math\'ematiques d'Orsay, Orsay Cedex, F-91405, France}
\email{simon.dauguet@math.u-psud.fr}

\date{\today}

\maketitle

\begin{abstract}
Dans cet article, on \'etend un crit\`ere d'ind\'ependance 
lin\'eaire d\^u \`a Fischler, qui est une g\'en\'eralisation 
quantitative du crit\`ere de Nesterenko, en affaiblissant 
fortement les hypoth\`eses sur les diviseurs des coefficients 
des formes lin\'eaires et en autorisant (dans une certaine mesure) ces formes \`a ne plus tendre vers 0. Ce nouveau 
crit\`ere est ensuite formul\'e dans un esprit plus \`a la Siegel 
en faisant intervenir une relation de r\'ecurrence v\'erifi\'ee 
par la suite de formes lin\'eaires. On en d\'emontre \'egalement une 
version plus g\'en\'erale, en termes de corps convexes et de 
r\'eseaux de $\R^n$.
\end{abstract}

\bigskip



\section{Introduction}


En 2000, Rivoal \cite{FctZetaInfValeursIrrationnelles} et Ball-Rivoal \cite{IrrationaliteInfValeursZeta} ont d\'emontr\'e qu'il existe une infinit\'e
d'entiers impairs en lesquels la fonction z\^eta de Riemann prend des valeurs irrationnelles. \`a la suite de ce th\'eor\`eme, une recherche du prochain $\zeta(2k+1)$
irrationnel
apr\`es la constante d'Ap\'ery $\zeta(3)$ a \'et\'e lanc\'ee. Apr\`es un premier r\'esultat de Rivoal \cite{Rivoal21}, Zudilin \cite{Irrationalite4} a d\'emontr\'e qu'au moins un des quatre nombres
$\zeta(5),\zeta(7),\zeta(9),\zeta(11)$ est
irrationnel.
Il a \'egalement montr\'e avec Fischler \cite{MathAnnalen} qu'il existe un entier impair $j\leq 139$ tel que $1,\zeta(3),\zeta(j)$ sont lin\'eairement ind\'ependants
sur $\Q$ (ce qui raffine des r\'esultats ant\'erieurs de Ball-Rivoal \cite{IrrationaliteInfValeursZeta} et Zudilin \cite{Irr145}). 

Dans les preuves de ces \'enonc\'es (et notamment du r\'esultat de Rivoal et Ball-Rivoal), le crit\`ere d'ind\'ependance lin\'eaire de Nesterenko \cite{Nesterenko}
joue un r\^ole d\'eterminant. 

\begin{theo}[Crit\`ere d'ind\'ependance lin\'eaire de Nesterenko] \label{Critere Nesterenko initial}
Soit $\xi_1,\,\dots,\,\xi_{p-1}$ des r\'eels, $p\geq 2$. Soient $0<\alpha<1$ et $\beta>1$. On consid\`ere $p$ suites d'entiers
$(\ell_{1,n})_n,\,\dots,\,(\ell_{p,n})_n$ telles que :
\begin{enumerate}
\item $\lim_{n\to\infty} \l|\sum_{i=1}^r \ell_{i,n}\xi_i+\ell_{p,n}\r|^{1/n}=\alpha$, 
\item pour tout $i\in\{1,\,\dots,\,p\}$, $\limsup_{n\to\infty} |\ell_{i,n}|^{1/n}\leq\beta$.
\end{enumerate}

Alors on a :

$$\dim_\Q\Vect_{\Q}(1,\xi_1,\,\dots,\,\xi_{p-1}) \geq 1-\frac{\log \alpha}{\log \beta}\,.$$
\end{theo}

Il fut g\'en\'eralis\'e \`a plusieurs reprises. Nesterenko lui-m\^eme \cite{Nesterenkopadic} en donne une version
$p$-adique. Chantanasiri \cite{Chantanasiri} a donn\'e une g\'en\'eralisation \`a $\C_p$. Töpfer \cite{Topfer1} a utilis\'e la m\'ethode de Nesterenko pour
d\'evelopper un
crit\`ere d'ind\'ependance alg\'ebrique et l'utilise \cite{Topfer2} pour d\'emontrer des r\'esultats d'ind\'ependance alg\'ebrique de valeurs de fonctions de Mahler. Bedulev
\cite{Bedulev} donne un crit\`ere pour des formes lin\'eaires \`a coefficients dans un corps de nombres. 

Dans \cite{MathAnnalen}, Fischler et Zudilin ont g\'en\'eralis\'e le crit\`ere de Nesterenko pour exploiter la pr\'esence de diviseurs des coefficients des formes
lin\'eaires. Fischler d\'emontre aussi dans \cite{NesterenkoOscillate} une version du crit\`ere de Nesterenko pour des formes lin\'eaires oscillantes. Et dans
\cite{Nesterenko4Vectors}, il en donne une version dans le cas d'une suite de formes lin\'eaires petites en plusieurs vecteurs. Il d\'emontre notamment le r\'esultat suivant :

\begin{theo}[Fischler, 2011]\label{Nesterenko avec diviseurs Fischler}
Soit $\xi_1,\dots,\xi_{p-1}\in\R$ avec $p\geq 2$. 

Soit $\tau_1,\dots,\tau_{p-1}>0$ des r\'eels deux \`a deux distincts et $\gamma_1,\dots,\gamma_{p}\geq 0$. Soit $(Q_n)_{n\in\N}$ une suite strictement croissante d'entiers
naturels telle que
$$
Q_{n+1}=Q_n^{1+o(1)}
$$
quand $n\to+\infty$. 
Pour tout $n\in\N$ et tout $i\in\{1,\dots,p\}$, soit $\ell_{i,n}\in\Z$ et $\delta_{i,n}$ un diviseur positif de $\ell_{i,n}$
tels que : 

\begin{enumerate}[(a)]
\item $\delta_{i,n}$ divise $\delta_{i+1,n}$ pour tous $n\in\N$ et tout $i\in\{1,\dots,p-1\}$.
\item $\frac{\delta_{j,n}}{\delta_{i,n}}$ divise $\frac{\delta_{j,n+1}}{\delta_{i,n+1}}$ pour tout $n\in\N$ et tout $0\leq
i<j\leq p$, avec $\delta_{0,n}=1$. 
\item pour tout $i\in\{1,\dots, p\}$, $\delta_{i,n}=Q_n^{\gamma_i+o(1)}$. 
\item $|\ell_{p,n}\xi_j-\ell_{j,n}|=Q_n^{-\tau_j+o(1)}\text{ pour tout }j\in\{1,\dots,p-1\}$.
\item $\max_{1\leq i\leq p} |\ell_{i,n}|\leq Q_n^{1+o(1)}$.
\end{enumerate}

Soient $\ve>0$, $Q$ suffisamment grand par rapport \`a $\ve$, et $(a_1,\dots,a_{p})\in\Q^{p}\setminus\{\mathbf0\}$ tels que pour tout
$i\in\{1,\dots,p\}$, $\delta_{i,\Phi(Q)}a_i\in\Z$ et pour tout $i\in\{1,\dots,p-1\}$, $|a_i|\leq Q^{\tau_i-\ve}$, o\`u $\Phi(Q)=\max\{k\in\N,Q_k\leq Q\}$. Alors
on a : 
$$
|a_1\xi_1+\dots+a_{p-1}\xi_{p-1}+a_p|>Q^{-1-\ve}\,.
$$

En particulier, les nombres $1,\xi_1,\dots,\xi_{p-1}$ sont lin\'eairement ind\'ependants sur $\Q$. 
\end{theo}

Cet \'enonc\'e est le Corollary 3 de \cite{Nesterenko4Vectors} dans le cas o\`u $\delta_{i,n}=1$ pour tous $i,n$. Le cas g\'en\'eral
se d\'eduit facilement du Theorem 3 de \cite{Nesterenko4Vectors} 
en prenant $k+1=p$, $\omega_j=\vphi_j=0$, $v_i=(0,\dots,0,1,0,\dots,0,-\xi_i)$ avec le 1 en $i$-\`eme position pour $1\leq i\leq p-1$ et
$v_p=(\xi_1,\dots,\xi_{p-1},1)$.

Dans cet article, on g\'en\'eralise le th\'eor\`eme \ref{Nesterenko avec diviseurs Fischler} en d\'emontrant notamment le r\'esultat suivant : 

\begin{theo}\label{Nesterenko avec diviseurs}
Soient $\xi_1,\dots,\xi_{p-1}$ des r\'eels. Soient $\tau_{1},\dots,\tau_{p-1}>-1$ des r\'eels 2 \`a 2 distincts et
$(Q_n)_{n\in\N}$ une suite d'entiers
strictement croissante telle que $Q_{n+1}=Q_n^{1+o(1)}$. Pour tout $i\in\{1,\dots,p\}$ et pour tout $n\in\N$, soit
$\delta_{i,n}\in\N^*$ 
tel que $\delta_{i,n}|\delta_{i,n+1}$. 
Soient maintenant des entiers $\ell_{i,n}$ pour $i=1,\dots,p$ et pour tout $n\in\N$, tels que : 
\begin{enumerate}[(i)]
  \item $\forall i\in\{1,\dots, p\}$, $\delta_{i,n}|\ell_{i,n}$ ;
  \item $\forall i\in\{1,\dots, p-1\}$, $|\ell_{i,n}-\ell_{p,n}\xi_i|=Q_n^{-\tau_i+o(1)}$ ; 
  \item $|\ell_{p,n}|=Q_n^{1+o(1)}$. 
\end{enumerate}

Soient $\ve>0$, $Q$ suffisamment grand par rapport \`a $\ve$ et $(a_1,\dots,a_p)\in\Q^p\setminus\{\textbf{0}\}$ tel que pour tout $i\in\{1,\dots,p\}$, $\delta_{i,\Phi(Q)}a_i\in\Z$ et pour tout $i\in\{1,\dots,p-1\}$, $|a_i|\leq Q^{\tau_i-\ve}$, o\`u $\Phi(Q)=\max\{k\in\N, Q_k\leq Q\}$. Alors on a :
$$
|a_1\xi_1+\dots+a_{p-1}\xi_{p-1}+a_p|>Q^{-1-\ve}\,.
$$
\end{theo}

La diff\'erence majeure avec le th\'eor\`eme \ref{Nesterenko avec diviseurs Fischler} de Fischler est certainement le fait qu'ici, les $\tau_i$ peuvent \^etre
n\'egatifs, tant qu'ils restent $>-1$. Auparavant, les formes lin\'eaires en $1$ et $\xi_i$ devaient tendre vers 0, alors que maintenant, elles peuvent tendre vers
$+\infty$. En revanche, on perd l'ind\'ependance lin\'eaire des $\xi_i$ qui provenait justement du fait d'avoir les $\tau_i>0$. D'autre part, les conditions de
divisibilit\'e sur les diviseurs $\delta_{i,n}$ ont \'et\'e fortement affaiblies. On ne demande m\^eme plus de comportement asymptotique sp\'ecifique pour ces suites de
diviseurs. En revanche, les comportements asymptotiques des suites de formes lin\'eaires ont \'et\'e renforc\'es : dans le th\'eor\`eme \ref{Nesterenko avec diviseurs
Fischler} de Fischler, on ne demande qu'une majoration alors que dans le th\'eor\`eme \ref{Nesterenko avec diviseurs} on exige que les coefficients soient
exactement de taille $Q_n^{1+o(1)}$.

\bigskip

Le th\'eor\`eme \ref{Nesterenko avec diviseurs} est d\'emontr\'e au \S~\ref{Applications} ci-dessous, comme cas particulier d'un r\'esultat plus g\'en\'eral qui fait l'objet du \S~\ref{criteres generaux} et qui concerne une base et une suite de r\'eseaux quelconques. Par ailleurs, dans le cas g\'en\'eral comme dans le cas particulier, on d\'emontre des variantes ``\`a la Siegel'' de ces crit\`eres. Enfin, au \S~\ref{Optimalite}, on \'etudie l'optimalit\'e de ces crit\`eres et on d\'emontre une r\'eciproque partielle. 

\bigskip

Ces crit\`eres ont \'et\'e appliqu\'es dans \cite{DauZud}, en collaboration avec Zudilin : en construisant des approximations simultan\'ees de
$\zeta(2)$
et $\zeta(3)$ \`a l'aide d'outils hyperg\'eom\'etriques, ils permettent d'obtenir une mesure d'ind\'ependance lin\'eaire restreinte de 1, $\zeta(1)$ et $\zeta(3)$ (sachant que l'ind\'ependance lin\'eaire sur $\Q$ de ces trois nombres reste conjecturale).  

Certains d\'etails suppl\'ementaires pourront \^etre obtenus dans la th\`ese de l'auteur \cite{these}, notamment sur l'application des r\'esultats montr\'es ici \`a des
formes lin\'eaires en 1, $\zeta(2)$ et $\zeta(3)$. 

\bigskip

Pour le reste de cet article, $p$ sera un entier $\geq2$.


\section{Crit\`eres quantitatifs g\'en\'eraux}\label{criteres generaux}


Le th\'eor\`eme suivant g\'en\'eralise le th\'eor\`eme \ref{Nesterenko avec diviseurs} de l'introduction (qui est obtenu dans le cas
particulier o\`u $e_i=(0,\dots,0,1,0,\dots,-\xi_i)$ avec le $1$ en $i$-\`eme position pour $1\leq i\leq p-1$ et $e_p=(\xi_1,\dots,\xi_{p-1},1)$, avec le r\'eseau $\Lambda_n$ de la forme $\bigoplus_{i=1}^p \delta_{i,n}\Z$). 

On note $(\R^p)^*$ le dual de $\R^p$, sur lequel on fixe une norme $\|.\|$.

\begin{theo}\label{prop a la Nesterenko general}

Soit $(e_1,\dots,e_p)$ une base de $\R^p$.

Soient $\tau_{1},\dots,\tau_{p-1}$ des r\'eels $>-1$ et 2 \`a 2 distincts. Soit $(Q_n)_{n\in\N}$ une suite d'entiers strictement
croissante telle que
$Q_{n+1}=Q_n^{1+o(1)}$.
Soit $(\Lambda_n)_{n\in\N}$
une suite
de r\'eseaux de $(\R^p)^*$ telle que $\Lambda_{n+1}\subset\Lambda_n$. On consid\`ere une suite $(L_n)_{n\in\N}$ de formes
lin\'eaires telle que
\begin{enumerate}
\item \label{1} $L_n\in\Lambda_n$ pour tout $n\in\N$,
\item \label{2} $|L_n(e_i)| = Q_n^{-\tau_i+o(1)}$ pour $1\leq i\leq p-1$, 
\item \label{3} $\|L_n\| = Q_n^{1+o(1)}$. 
\end{enumerate}
Alors pour tout $\ve>0$ et pour tout $Q$ assez grand en termes de $\ve$, on a
$$
\Lambda_{\Phi(Q)}^\perp\cap\mathcal{C}_{Q,\ve}=\{\textbf{0}\}\,,
$$
o\`u 
$\Phi(Q)=\max\{k\in\N,Q_k\leq Q\}$, $\Lambda_n^\perp$ est le r\'eseau dual de $\Lambda_n$ et  

$$
\mathcal{C}_{Q,\ve}=\l\{\sum_{i=1}^p \lambda_ie_i, |\lambda_i|<Q^{\tau_i-\ve}\text{ pour }1\leq i\leq p-1\text{ et
}|\lambda_p|<Q^{-1-\ve}\r\}\,.
$$
\end{theo}

Par d\'efinition, $\Lambda_n^\perp$ est l'ensemble des vecteurs $x\in\R^p$ tels que $L(x)\in\Z$ pour toute forme lin\'eaire $L\in\Lambda_n$ (voir par exemple
\cite{Cassels} chapitre I.5 page 23). Dans le cas particulier de la section \ref{Applications} ci-dessous, $\Lambda_n=\bigoplus_{i=1}^p \delta_{i,n}\Z$ est
l'ensemble des formes lin\'eaires $L=\ell_1X_1+\dots+\ell_pX_p$ telles que $\delta_{i,n}\big|\ell_i$ pour tout $1\leq i\leq p$. On voit alors que
$\Lambda_n^\perp$ est l'ensemble des $(x_1,\dots,x_p)\in\Q^p$ tels que $\delta_{i,n}x_i\in\Z$ pour tout $i\in\{1,\dots,p\}$.

\bigskip

La preuve \cite{Nesterenko4Vectors} du th\'eor\`eme \ref{Nesterenko avec diviseurs Fischler} repose sur le lemme suivant (Lemma 3
de \cite{Nesterenko4Vectors})~:

\begin{lem}[ Fischler, 2011] \label{lemme Fischler Nesterenko4Vectors}
Soit $M\in\MM_p(\R)$ une matrice carr\'ee de taille $p$, dont les coefficients $m_{i,j}$ sont non nuls et v\'erifient 
\begin{equation}
|m_{i',j}m_{i,j'}|\leq\frac{1}{(p+1)!}|m_{i,j}m_{i',j'}|\text{ pour tous }i,j,i',j'\text{ tels que }i<i'\text{ et }j<j'\label{mm'}\,.
\end{equation}
Alors $M$ est inversible.
\end{lem}

Ce lemme est d\'emontr\'e et utilis\'e dans \cite{Nesterenko4Vectors} (avec une majoration des coeficients de $M^{-1}$ qui est inutile ici) sous l'hypoth\`ese que
les coefficients $m_{i,j}$ sont strictement positifs. Mais la preuve de \cite{Nesterenko4Vectors} fonctionne plus g\'en\'eralement d\`es que $m_{i,j}\not=0$, car elle
utilise seulement les majorations \eqref{mm'} sur les valeurs absolues des coefficients.

\begin{proof}[D\'emonstration du th\'eor\`eme \ref{prop a la Nesterenko general}]
Tout d'abord, notons $|||.|||$ la norme sur $(\R^p)^*$ subordonn\'ee \`a la norme $N$ sur $\R^p$ d\'efinie par $N(\mu_1e_1+\dots+\mu_pe_p)=\max_{1\leq i\leq p}
|\mu_i|$. Pour tout $n\in\N$, il existe $x=\mu_1e_1+\dots+\mu_pe_p$ tel que $N(x)=1$ et 
$$
|||L_n|||=|L_n(x)|\leq \sum_{i=1}^p |\mu_i L_n(e_i)|\leq |L_n(e_p)|+\sum_{i=1}^{p-1}Q_n^{-\tau_i+o(1)}\,.
$$
Comme $\tau_i>-1$ pour tout $1\leq i \leq p-1$, l'hypoth\`ese \eqref{3} montre (puisque toutes les normes sont \'equivalentes sur $(\R^p)^*$) que $|L_n(e_p)|\geq
Q_n^{1+o(1)}$. Comme l'in\'egalit\'e dans l'autre sens d\'ecoule imm\'ediatement de l'hypoth\`ese \eqref{3}, on a \'egalit\'e. 

Soit $\ve>0$ et $Q$ suffisamment grand par rapport \`a $\ve$. Soit $\ve_1>0$ tel que 
\begin{equation} \label{condition ve_1 general}
\l\{\begin{matrix} 
\tau_i\l(1-(1+\ve_1)^{p-1}\r)<\frac{\ve}{2} & \text{ pour tout }i\in\{1,\dots,p-1\} \\
(\ve_1+1)^{p-1}-1<\frac{\ve}{2} 
\end{matrix}\r.\,.
\end{equation}
On remarque que seuls les $i\in\{1,\dots,p-1\}$ tels que $\tau_i<0$ fournissent une contrainte sur $\ve_1$. 

On d\'efinit $\varphi:\N^*\to\N^*$ par $\varphi(n)-1=\Phi(Q_n^{1+\ve_1})=\max\{k\in\N, Q_k\leq Q_n^{1+\ve_1}\}$. On en d\'eduit $\varphi(n)\geq n+1$, d'o\`u
$\lim_{n\to\infty}\varphi(n)=+\infty$. Ainsi, $Q_{\varphi(n)}=Q_{\varphi(n)-1}^{1+o(1)}$. Mais d'autre part,
la d\'efinition de
$\varphi$ nous donne 
$$
Q_{\vphi(n)}=Q_n^{1+\ve_1+o(1)}\,,
$$
de sorte que 
$$
Q_{\vphi_i(n)}=Q_{\vphi_{i-1}(n)}^{1+\ve_1+o(1)}=Q_n^{(1+\ve_1)^i+o(1)}\,,
$$
o\`u $\vphi_i$ est la $i$-\`eme it\'er\'ee de $\vphi$ (\emph{i.e.}, $\vphi_i=\underbrace{\vphi\circ\dots\circ\vphi}_{i\text{ fois}}$), avec $\varphi_0(n)=n$ et
$\varphi_1(n)=\varphi(n)$.

On choisit maintenant $n=\Phi(Q)$. Par d\'efinition de $\Phi(Q)$, on a donc $Q_n\leq Q<Q_{n+1}\leq Q_{\vphi(n)}$. On remarque alors que
$n\xrightarrow[Q\to\infty]{}+\infty$ et donc $o(1)$ est une suite tendant vers 0 quand $n$ ou $Q$ tendent vers $+\infty$ indiff\'eremment. On pourra donc utiliser $o(1)\xrightarrow[Q\to\infty]{}0$ ou $o(1)\xrightarrow[n\to\infty]{}0$ sans distinction ; on notera
simplement $o(1)$. De
plus, l'encadrement pr\'ec\'edent nous donne aussi $Q_n=Q^{1+o(1)}$. 

En posant $M_n=(L_{\vphi_{i-1}(n)}(e_j))_{1\leq i,j\leq p}\in\MM_p(\R)$, on a 
$$
\l|\frac{L_{\vphi_{i'-1}(n)}(e_j)L_{\vphi_{i-1}(n)}(e_{j'})}{L_{\vphi_{i-1}(n)}(e_j)L_{\vphi_{i'-1}(n)}(e_{j'})}\r| =
Q_n^{-(\tau_{i'}-\tau_i)((1+\ve_1)^{j'}-(1+\ve_1)^j)+o(1)}\,,
$$
et donc le lemme \ref{lemme Fischler Nesterenko4Vectors} s'applique pour $n$ suffisamment grand. Ainsi, la matrice $M_n$ est inversible et les formes lin\'eaires
$L_n,\dots,L_{\vphi_{p-1}(n)}$ sont lin\'eairement ind\'ependantes.

On consid\`ere 
$P\in\mathcal{C}_{Q,\ve}\cap\Lambda_{n}^\perp$ non nul. Si $Q$ est assez grand, $L_n,\dots,L_{\varphi_{p-1}(n)}$ sont lin\'eairement
ind\'ependants et forment donc une base de l'espace dual de $\R^p$ : il existe $k\in\{0,\dots,p-1\}$ tel que
$L_{\varphi_k(n)}(P)\not=0$. 

Par hypoth\`ese sur la suite de r\'eseaux, on a $\Lambda_n^{\perp}\subset\Lambda_{\varphi^k(n)}^\perp$. Donc, par d\'efinition du r\'eseau dual,
$L_{\varphi_k(n)}(P)\in\Z\setminus\{0\}$. 

Mais en notant $P=\sum_{i=1}^p \lambda_ie_i$ comme dans la d\'efinition de $\mathcal{C}_{Q,\ve}$, on a :

\begin{align*}
|L_{\varphi_k(n)}(P)| 
&\leq \sum_{i=1}^p |\lambda_i||L_{\varphi_k(n)}(e_i)| \\
&= \sum_{i=1}^{p-1} Q^{\tau_i-\ve}\l(Q_n^{(1+\ve_1)^k+o(1)}\r)^{-\tau_i+o(1)} + Q^{-1-\ve}\l(Q_n^{(1+\ve_1)^k+o(1)}\r)^{1+o(1)} \\
&\leq \sum_{i=1}^{p-1} Q^{\tau_i(1-(1+\ve_1)^k)-\ve/2} + Q^{(1+\ve_1)^k-1-\ve/2}\,.
\end{align*}

Mais par le choix \eqref{condition ve_1 general} de $\ve_1$, chacun des exposants est $<0$. Donc, en prenant $Q$ assez grand, on aura $|L_{\varphi_k(n)}(P)|<1$,
ce qui aboutit \`a une contradiction. 
\end{proof}


\bigskip 

Le crit\`ere d'ind\'ependance lin\'eaire de Nesterenko~\cite{Nesterenko} admet un analogue ``\`a la Siegel'' dont la preuve repose simplement
sur un calcul de d\'eterminant (voir par exemple la proposition 1 de~\cite{Nesterenko4Vectors} ou la proposition 4.1 de Marcovecchio~\cite{Marcovecchio}). Dans le
cas
dual o\`u on part d'approximations simultan\'ees de $\xi_1,\dots,\xi_{p-1}$, on peut de m\^eme obtenir un analogue du th\'eor\`eme \ref{Nesterenko avec diviseurs
Fischler} ;
il s'agit (en l'absence de diviseurs) d'un cas particulier de la proposition 1 de~\cite{Nesterenko4Vectors}. Pour exploiter la pr\'esence de diviseurs
$\delta_{i,n}$, on peut appliquer le r\'esultat suivant ``\`a la Siegel'', dans lequel les conditions sur les tailles des objets sont affaiblies. Le th\'eor\`eme
n\'ecessite en contrepartie une quantit\'e de formes lin\'eaires cons\'ecutives et lin\'eairement ind\'ependantes ainsi qu'une relation de r\'ecurrence v\'erifi\'ee par la suite
de formes lin\'eaires. En effet, Siegel a \'et\'e le premier \`a comprendre dans son article \cite{Siegel}\footnote{On pourra trouver deux versions r\'esum\'ees de cet
article
dans \cite{NumberTheory} Chapitre 2, §1, pages 81-82 et Chapitre 5, §2, pages 215-216.} que l'on pouvait d\'emontrer l'ind\'ependance lin\'eaire de nombres
$\xi_1,\dots,\xi_{p-1}$ en construisant un syst\`eme de formes lin\'eaires lin\'eairement ind\'ependantes \`a coefficients entiers, petites en $(\xi_1,\dots,\xi_{p-1})$.

\begin{theo}\label{prop a la Siegel general}
Soit $(e_1,\dots,e_p)$ une base de $\R^p$.

Soient $\tau_{1},\dots,\tau_{p-1}$ des r\'eels $>-1$ et $(Q_n)_{n\in\N}$ une suite d'entiers strictement croissante telle que $Q_{n+1}=Q_n^{1+o(1)}$.
Soit $(\Lambda_n)_{n\in\N}$
une suite
de r\'eseaux de $(\R^p)^*$ telle que $\Lambda_{n+1}\subset\Lambda_n$. On consid\`ere une suite $(L_n)_{n\in\N}$ de formes lin\'eaires telle que 
\begin{enumerate}[(i)]
\item pour tout $n\in\N$, $L_n\in\Lambda_n$ ; 
\item $|L_n(e_i)|\leq Q_n^{-\tau_i+o(1)}$ pour $1\leq i\leq p-1$ ; 
\item $\|L_n\|\leq Q_n^{1+o(1)}$ ;
\item la suite de formes
lin\'eaires v\'erifie une relation de r\'ecurrence : il existe $n_1\geq 0$ tel que pour tout $n\geq n_1$,
il existe des nombres r\'eels $\alpha_0(n),\dots, \alpha_{p-1}(n)$ tels que $L_{n+p}=\sum_{i=0}^{p-1} \alpha_i(n)L_{n+i}$ avec $\alpha_0(n)\not=0$. 
\end{enumerate}

On suppose enfin qu'il existe $n_2\geq n_1$ tel que $L_{n_2},\dots,L_{n_2+p-1}$ soient lin\'eairement ind\'ependantes.

Alors pour tout $\ve>0$, pour tout $Q$ assez grand en termes de $\ve$, on a 
$$
\Lambda_{\Phi(Q)}^\perp\cap\mathcal{C}_{Q,\ve}=\{\textbf{0}\}\,,
$$
o\`u $
\Phi(Q)=\max\{k\in\N, Q_k\leq Q\}$ et
$$
\mathcal{C}_{Q,\ve}=\l\{\sum_{i=1}^p \lambda_ie_i, |\lambda_i|\leq Q^{\tau_i-\ve}\text{ pour }1\leq i\leq p-1\text{, et }|\lambda_p|\leq Q^{-1-\ve}\r\}\,.
$$
\end{theo}

Notons qu'\`a
la place de la relation de r\'ecurrence et de l'existence de $p$ formes lin\'eaires cons\'ecutives lin\'eairement ind\'ependantes, on peut supposer que parmi $k$ formes
lin\'eaires cons\'ecutives
$L_n,\dots,L_{n+k-1}$, il y en a
toujours $p$ qui sont lin\'eairement ind\'ependantes (o\`u $k\geq p$ est un entier fix\'e). En effet, on peut alors noter $\vphi$ l'extraction telle que
$\vphi(n+1)=\min\l\{m>\vphi(n),L_m\notin\Vect(L_{\vphi(n-p+2)},\dots,L_{\vphi(n)})\r\}$, et on a
$\vphi(n+1)-\vphi(n)\leq k-1$ qui est une constante, de sorte que $Q_{\vphi(n+1)}=Q_{\vphi(n)}^{1+o(1)}$. Comme $L_{\vphi(n)},\dots,L_{\vphi(n+p-1)}$ est une
base de $(\R^p)^*$, on peut \'ecrire $L_{\vphi(n+p)}=\sum_{j=0}^{p-1} \alpha_j(n)L_{\vphi(n+j)}$ et on a n\'ecessairement $\alpha_0(n)\not=0$, car sinon
$L_{\vphi(n+p)}\in\Vect(L_{\vphi(n+1)},\dots,L_{\vphi(n+p-1)})$ ce qui n'est pas possible par choix de
$\vphi$. Et le th\'eor\`eme s'applique donc \`a la suite extraite grâce \`a $\vphi$.

Quand on suppose qu'il y a une relation de r\'ecurrence, son ordre doit \^etre n\'ecessairement $\geq p$ car l'hypoth\`ese de $p$ formes lin\'eaires lin\'eairement ind\'ependantes
est primordiale. Cette
hypoth\`ese nous permet d'assurer que $L_{n+k}(P)$ n'est pas nulle pour un certain $k\in\{0,\dots, p-1\}$ d\`es que $P\not=\textbf{0}$, point essentiel de la preuve.

\begin{proof}
Comme au d\'ebut de la preuve de la proposition \ref{prop a la Nesterenko general}, on remarque que $|L_n(e_p)|\leq Q_n^{1+o(1)}$. 

La matrice
$$
\Delta_{n}=\begin{pmatrix} L_{n}(e_1) & \dots & L_{n+p-1}(e_1) \\ \vdots & \ddots & \vdots \\ L_{n}(e_{p}) & \dots & L_{n+p-1}(e_{p})
\end{pmatrix}
$$
est de d\'eterminant non nul pour $n=n_2$. D'autre part, on remarque que
\begin{align*}
\rg(\Delta_{n+1}) 
&= \rg(\Delta_{n})
\end{align*}
car $\alpha_0(n)\not=0$ donc pour tout $n\geq n_2$, $\rg(\Delta_n)=\rg(\Delta_{n_2})=p$. Donc toutes les matrices $\Delta_n$ sont de
d\'eterminant non nul. 

De m\^eme que dans la d\'emonstration du th\'eor\`eme \ref{prop a la Nesterenko general}, on a $n\xrightarrow[Q\to+\infty]{} +\infty$
si $n=\Phi(Q)$, de sorte que
$Q_n=Q^{1+o(1)}$ et les suites $o(1)$ sont des suites tendant vers 0 quand $n$ ou $Q$ tendent vers $+\infty$ indiff\'eremment. 

Soit $\ve>0$ et $Q$ assez grand en termes de $\ve$ pour avoir $n=\Phi(Q)\geq n_1,n_2$. Soit
$P\in\Lambda_{\Phi(Q)}^\perp\cap\mathcal{C}_{Q,\ve}\setminus\{\mathbf0\}$. Comme $\det(\Delta_n)\not=0$, on sait que toute combinaison
lin\'eaire non triviale
des lignes de la matrice $\Delta_n$ est non nulle. En particulier, il existe $k\in\{0,\dots,p-1\}$ tel que $L_{n+k}(P)\not=0$. En outre, $L_{n+k}(P)\in\Z$ car $P\in\Lambda_n^\perp\subset\Lambda_{n+k}^\perp$.

La fin de la preuve est alors identique \`a celle du th\'eor\`eme \ref{prop a la Nesterenko general}.
\end{proof}


\section{Applications \`a une base particuli\`ere de $\R^p$ et \`a des r\'eseaux particuliers}\label{Applications}


Dans toute cette section, on se place dans le cadre particulier d'une base de la forme 
$e_i=(0,\dots,0,1,0,\dots,0,-\xi_i)$ avec le 1 en $i$-\`eme position pour $1\leq i\leq p-1$ et $e_p=(\xi_1,\dots,\xi_{p-1},1)$ et de r\'eseaux
$\Lambda_n=\bigoplus_{i=1}^p \delta_{i,n}\Z$. Ce cas particulier a d\'ej\`a \'et\'e mentionn\'e au d\'ebut de la section \ref{criteres generaux}.

Le th\'eor\`eme \ref{Nesterenko avec diviseurs} de l'introduction est en fait un corollaire du th\'eor\`eme \ref{prop a la Nesterenko general} (voir la d\'emonstration
ci-dessous). Son int\'er\^et est qu'il est plus conforme \`a ce qu'on a en pratique. 

Il raffine le th\'eor\`eme \ref{Nesterenko avec diviseurs Fischler} de \cite{Nesterenko4Vectors} cit\'e dans l'introduction : la principale nouveaut\'e est que les
formes lin\'eaires
$\ell_{i,n}-\xi_i\ell_{p,n}$ ne tendent plus forc\'ement vers 0, puisque les $\tau_i$ peuvent \^etre n\'egatifs. Par ailleurs, les contraintes sur les diviseurs
$\delta_{i,n}$ sont affaiblies. La seule contrepartie est qu'on demande aux coefficients $\ell_{i,n}$ d'\^etre, en valeur absolue, \'egaux \`a $Q_n^{1+o(1)}$ (il
suffit de faire cette hypoth\`ese pour $i=p$, puisque $\tau_i>-1$ pour tout $i\leq p-1$). En pratique, cela ne pose pas de probl\`eme (quitte \`a diminuer $Q_n$ si
n\'ecessaire, ce qui am\'eliore la conclusion du th\'eor\`eme), puisqu'on a $|\ell_{p,n}|<|\ell_{p,n+1}|=|\ell_{p,n}|^{1+o(1)}$ dans les
applications.

Pour d\'emontrer le th\'eor\`eme \ref{Nesterenko avec diviseurs Fischler}, Fischler combine le lemme \ref{lemme Fischler Nesterenko4Vectors} (avec $p-1$ au lieu de
$p$) avec
le premier th\'eor\`eme de Minkowski
sur les corps convexes. Pour d\'emontrer le th\'eor\`eme \ref{Nesterenko avec diviseurs}, on applique ce lemme \`a une matrice carr\'ee de taille $p$ (en utilisant
notamment l'hypoth\`ese (iii)) et on obtient ainsi $p$ vecteurs lin\'eairement ind\'ependants qui fournissent des approximations simultan\'ees de
$\xi_1,\dots,\xi_{p-1}$
par des nombres rationnels ayant le m\^eme d\'enominateur. Cela permet de conclure la preuve par un argument \`a la Siegel, sans utiliser la g\'eom\'etrie des nombres.
On \'evite ainsi les hypoth\`eses sur les diviseurs, pr\'esentes dans le th\'eor\`eme \ref{Nesterenko avec diviseurs Fischler}, et on a pas
besoin de supposer $\tau_1,\dots,\tau_{p-1}>0$. 

Une preuve directe du th\'eor\`eme \ref{Nesterenko avec diviseurs} figure dans le paragraphe 4 de \cite{these}. Dans cet article, on d\'eduit ce r\'esultat du th\'eor\`eme
\ref{prop a la Nesterenko general} d\'emontr\'e au \S~\ref{criteres generaux}. 

On remarque aussi que sous l'hypoth\`ese suppl\'ementaire $\delta_{i,n}=Q_n^{\gamma_i+o(1)}$ (hypoth\`ese qui sera utile au \S~\ref{Optimalite}), on a
n\'ecessairement $a_i=0$ pour tous les $i$ tels que $\tau_i+\gamma_i<0$. Le lecteur obtiendra plus de d\'etails dans
la d\'emonstration de la propositions \ref{construction ell et a}. 

Outre les formes lin\'eaires qui peuvent tendrent vers $+\infty$ au lieu de seulement 0, les diff\'erences entre le th\'eor\`eme \ref{Nesterenko avec diviseurs} et le th\'eor\`eme \ref{Nesterenko avec diviseurs Fischler} sont les suivantes~:
\begin{itemize}
\item Dans le th\'eor\`eme \ref{Nesterenko avec diviseurs}, on ne d\'emontre plus que
$1,\xi_1,\dots,\xi_{p-1}$ sont lin\'eairement ind\'ependants sur $\Q$. 
C'est une contrepartie de la
remarque pr\'ec\'edente selon laquelle les formes lin\'eaires ne
tendent pas forc\'ement vers 0. La conclusion du th\'eor\`eme \ref{Nesterenko avec diviseurs} n'est donc que quantitative. 

\item Si $\delta_{i,n}=Q_n^{\gamma_i+o(1)}$ quand $n\to+\infty$ avec $1\leq i\leq p-1$ et $\gamma_i\geq 0$ tel que
$\tau_i\leq -\gamma_i$, alors
l'entier
$\delta_{i,\Phi(Q)}a_i$ est major\'e par $Q^{-\ve+o(1)}$, donc $a_i=0$ d\`es que $Q$ est assez grand. Dans ce cas, le th\'eor\`eme \ref{Nesterenko avec diviseurs} ne
concerne pas vraiment $\xi_i$. Lorsque $\gamma_i=0$ pour tout $i$ (notamment si $\delta_{i,n}=1$), la situation est identique \`a celle du
th\'eor\`eme \ref{Nesterenko avec diviseurs Fischler} : la conclusion concerne uniquement les $\xi_i$ tels que $\tau_i>0$. 

\item La condition (a) du th\'eor\`eme \ref{Nesterenko avec diviseurs Fischler} sur les diviseurs a
disparu et la condition (b) sur la divisibilit\'e des quotients successifs a \'et\'e restreinte \`a $i=0$. 

\item La condition asymptotique (c) sur les suites de diviseurs $(\delta_{i,n})_{n\geq 0}$ a disparu. En revanche, la
condition (e) sur la taille des coefficients
de la suite des formes lin\'eaires a
\'et\'e renforc\'ee~: dans le th\'eor\`eme \ref{Nesterenko avec diviseurs Fischler}, seule l'in\'egalit\'e $\max_{0\leq i\leq p}
|\ell_{i,n}|\leq
Q_n^{1+o(1)}$
lorsque $n\to+\infty$ est exig\'ee.
Mais dans le th\'eor\`eme \ref{Nesterenko avec diviseurs}, on exige que la taille soit exactement de $Q_n^{1+o(1)}$. En d'autres termes, on change cette
hypoth\`ese en la condition plus forte $|\ell_{p,n}|=Q_n^{1+o(1)}$ quand $n\to+\infty$. Cette condition est cruciale dans la preuve et elle est souvent v\'erifi\'ee
en pratique.
\end{itemize}

\begin{proof}[D\'emonstration du th\'eor\`eme \ref{Nesterenko avec diviseurs}]
Comme indiqu\'e dans le d\'ebut de cette partie, on pose \break $e_i=(0,\dots,0,1,0,\dots,0,-\xi_i)$ avec le 1 en $i$-\`eme position, pour $1\leq i\leq p-1$, et
$e_p=(\xi_1,\dots,\xi_{p-1},1)$ ; alors $(e_1,\dots,e_p)$ forme une base de $\R^p$. 

On pose maintenant $\ve'=\ve/2$ et 
\begin{equation}\label{Q'}
Q'=\l(\l(1+\sum_{i=1}^{p-1}\xi_i^2\r)Q^{1+\ve}\r)^{\frac{1}{1+\ve'}}\,.
\end{equation}
$Q'$ peut donc \^etre pris aussi grand que n\'ecessaire par rapport \`a $\ve'$, en choisissant $Q$ assez grand.

On pose, pour tout $n\in\N$, $\Lambda_n=\bigoplus_{i=1}^p \delta_{i,n}\Z$. La relation $\delta_{i,n}\big|\delta_{i,n+1}$ nous donne
$\Lambda_{n+1}\subset\Lambda_n$. On pose aussi $L_n=(\ell_{i,n},\dots,\ell_{p,n})$. On a donc $L_n\in\Lambda_n$ pour tout $n\in\N$. On a \'egalement, pour tout
$i\in\{1,\dots,p-1\}$ et tout $n\in\N$, $|L_n(e_i)|=|\ell_{i,n}-\ell_{p,n}\xi_i|=Q_n^{-\tau_i+o(1)}$. L'in\'egalit\'e triangulaire $|\ell_{i,n}|\leq
|\ell_{i,n}-\xi_i\ell_{p,n}|+|\xi_i||\ell_{p,n}|$ nous donne $|\ell_{i,n}|\leq Q_n^{1+o(1)}$ pour tout $i\in\{1,\dots,p-1\}$ car $\tau_i>-1$ de sorte que
$\|L_n\|_\infty\leq Q_n^{1+o(1)}$, o\`u $\|x\|_\infty=\max_{1\leq i\leq p} |x_i|$. Finalement, on obtient grâce \`a l'hypoth\`ese (iii), puisque toutes les normes sont \'equivalentes sur
$(\R^p)^*$, $\|L_n\|=Q_n^{1+o(1)}$. 

On peut donc appliquer le th\'eor\`eme \ref{prop a la Nesterenko general} pour obtenir, si $Q$ est assez grand,
\begin{equation}\label{enonce 1}
\Lambda_{\Phi(Q')}^\perp\cap\CC_{Q',\ve'}=\{\textbf{0}\}\,.
\end{equation}

On prend $(a_1,\dots,a_p)\in\Lambda_{\Phi(Q)}^\perp$ non nul, avec $|a_i|\leq
Q^{\tau_i-\ve}$ pour $1\leq i\leq p-1$ et
$|a_1\xi_1+\dots+a_{p-1}\xi_{p-1}+a_p|< Q^{-1-\ve}$. On va montrer qu'on aboutit \`a une contradiction. 

On pose, pour $i\in\{1,\dots,p-1\}$, 
$$
\lambda_i=a_i-\xi_i\frac{\sum_{k=1}^{p-1}a_k\xi_k+a_p}{1+\sum_{k=1}^{p-1}\xi_k^2}
$$
et 
$$
u=\frac{\sum_{1\leq k\leq p-1} a_k\xi_k+a_p}{1+\sum_{1\leq k\leq p-1}\xi_k^2}\,.
$$

On obtient alors les relations suivantes pour $i\in\{1,\dots,p-1\}$ : $a_i=\lambda_i+u\xi_i$ et
$a_p=u-\sum_{1\leq i\leq p-1}\lambda_i\xi_i$. On pose 
$P=\lambda_1e_1+\dots+\lambda_{p-1}e_{p-1}+ue_p$, si bien que $P=(a_1,\dots,a_{p-1},a_p)\in\Lambda_{\Phi(Q)}^\perp\setminus\{\textbf{0}\}$ par construction. 

On a les majorations suivantes :
\begin{align*}
|u| 
&\leq \frac{1}{1+\sum_{1\leq i\leq p-1}\xi_i^2} Q^{-1-\ve} =Q'^{-1-\ve'}
\end{align*}
et, pour $i\in\{1,\dots,p-1\}$ : 
\begin{align*}
|\lambda_i|
&\leq Q^{\tau_i-\ve} + \frac{|\xi_i|}{1+\sum_{k=1}^{p-1} \xi_k^2}Q^{-1-\ve} \\
&\leq \l(1+\frac{|\xi_i|}{1+\sum_{k=1}^{p-1} \xi_k^2}\r)Q^{\tau_i-\ve} &\text{car }\tau_i>-1 \\
&\leq \l(1+\sum_{k=1}^{p-1} \xi_k^2\r)^{\frac{\tau_i-\ve/2}{1+\ve}}Q^{\tau_i-\frac{\ve}{2}} & \text{si }Q\text{ est assez grand} \\
&= \l(\l(1+\sum_{k=1}^{p-1}\xi_k^2\r)Q^{1+\ve}\r)^{\frac{\tau_i-\ve/2}{1+\ve}} \\
&= Q'^{\frac{1+\ve'}{1+\ve}\l(\tau_i-\frac{\ve}{2}\r)}\leq Q'^{\tau_i-\ve'}
\end{align*}
par d\'efinition de $Q'$ et car $\frac{1+\ve'}{1+\ve}=\frac{1+\ve'}{1+2\ve'}<1$ et $\ve/2=\ve'$. 

Finalement, on a $P=\sum_{1\leq i\leq p-1}\lambda_ie_i+ue_p\not=\textbf{0}$ dans $\Lambda_{\Phi(Q)}^\perp$ avec $|\lambda_i|\leq Q'^{\tau_i-\ve'}$ et $|u|\leq
Q'^{-1-\ve'}$. Et comme $\frac{1+\ve}{1+\ve'}>1$, on a $Q'=\l(1+\sum_{i=1}^{p-1}\xi_i^2\r)^\frac{1}{1+\ve'}Q^{\frac{1+\ve}{1+\ve'}}
>Q$ si $Q$ est assez grand. Donc, $\Phi(Q')\geq\Phi(Q)$, et pour tout $i\in\{1,\dots,p-1\}$,
$\delta_{i,\Phi(Q)}$ divise $\delta_{i,\Phi(Q')}$. Ainsi $\Lambda_{\Phi(Q)}^\perp\subset\Lambda_{\Phi(Q')}^\perp$ d'o\`u $P\in\Lambda_{\Phi(Q')}^\perp$ et on a
une
contradiction avec \eqref{enonce 1}.
\end{proof}

En fait, l'implication d\'emontr\'ee ci-dessus entre la conclusion du th\'eor\`eme \ref{prop a la Nesterenko general} et celle du th\'eor\`eme \ref{Nesterenko
avec diviseurs} est une \'equivalence. Le lecteur pourra trouver la preuve pr\'ecise de cette \'equivalence dans la th\`ese de l'auteur \cite{these} (proposition 6.3 au paragraphe 6.2).

\bigskip


Comme dans la section \ref{criteres generaux}, on donne ici une version du th\'eor\`eme \ref{Nesterenko avec diviseurs} plus  ``\`a la Siegel'' dans le cadre
particulier de cette section.

\begin{theo}\label{Nesterenko a la Siegel}
Soient $\xi_1,\dots,\xi_{p-1}$ des r\'eels quelconques, et $\tau_1\,,\dots\,,\tau_{p-1}$ des r\'eels $>-1$. Soit $(Q_n)_{n\geq
0}$ une suite strictement
croissante d'entiers positifs telle que $Q_{n+1}=Q_n^{1+o(1)}$. On suppose \'egalement que pour tout $n\in\N$ et pour tout
$i\in\{1,\dots,p\}$, il existe
$\delta_{i,n}\in\N^*$ et $\ell_{i,n}\in\Z$ tels que $\delta_{i,n}$ divise $\delta_{i,n+1}$ et : 
\begin{enumerate}[(i)]
\item $\forall i\in\{1,\dots, p\}$, $\delta_{i,n}|\ell_{i,n}$ ;
\item $\forall i\in\{1,\dots, p-1\}$, $|\ell_{i,n}-\xi_i\ell_{p,n}|\leq Q_n^{-\tau_i+o(1)}$ ;
\item $|\ell_{p,n}|\leq Q_n^{1+o(1)}$ ;
\item il existe un entier $n_1$ tel que, pour tout $n\geq n_1$, il existe des r\'eels
$\alpha_0(n),\,\alpha_1(n),\,\dots,\,\alpha_{p-1}(n)$ avec
$\alpha_0(n)\not=0$, tels que pour tout $i\in\{1,\dots, p\}$ on ait :
  $$\ell_{i,n+p}=\sum_{j=0}^{p-1} \alpha_j(n)\ell_{i,n+j}\,.$$ 
\end{enumerate}
On suppose enfin que si l'on note $\Delta_n$ la matrice suivante de taille  $p\times p$ : 
$$
 \Delta_n=\begin{pmatrix} \ell_{1,n} & \cdots & \ell_{1,n+p-1} \\ \vdots & \ddots & \vdots \\ \ell_{p,n} & \cdots &
\ell_{p,n+p-1} \end{pmatrix}\ ,
$$
il existe un certain $n_2\geq n_1$ tel que $\det(\Delta_{n_2})\not=0$. 

Soit $\ve>0$, $Q>0$ suffisamment grand en fonction de $\ve$ et soit $(a_1,\dots,a_p)\in\Q^p\setminus\{\textbf{0}\}$, tel que
$\delta_{i,\Phi(Q)}a_i\in\Z$ pour
$1\leq i\leq p$ et $|a_i|\leq Q^{\tau_i-\ve}$ pour $1\leq i\leq p-1$, o\`u $\Phi(Q)=\max\{m\in\N, Q_m\leq Q\}$. Alors on a
$$
|a_1\xi_i+\dots+a_{p-1}\xi_{p-1}+a_p|\geq Q^{-1-\ve}\,.
$$
\end{theo}

De m\^eme que pour le th\'eor\`eme \ref{prop a la Siegel general}, il suffirait que parmi $k$ formes lin\'eaires cons\'ecutives $L_n,\dots,L_{n+k-1}$ avec $L_n=(\ell_{1,n},\dots,\ell_{p,n})$, il y en ait toujours $p$ qui soient lin\'eairement ind\'ependantes (o\`u $k\geq p$ est un entier fix\'e). 

La conclusion de ce crit\`ere ``\`a la Siegel'' est la m\^eme que celle du th\'eor\`eme \ref{Nesterenko avec diviseurs} (qui est davantage dans l'esprit du crit\`ere
de Nesterenko). Les diff\'erences entre ces deux \'enonc\'es sont les suivantes~:

\begin{itemize}
\item Dans tout crit\`ere \`a la Siegel, on a besoin d'une hypoth\`ese assurant l'ind\'ependance lin\'eaire de $p$ formes lin\'eaires. Dans le th\'eor\`eme \ref{Nesterenko a la
Siegel}, cette hypoth\`ese prend la forme d'un d\'eterminant non nul (celui de $\Delta_{n_2}$, avec une valeur fix\'ee de $n_2$) et d'une relation de
r\'ecurrence.
L'absence de cette hypoth\`ese est un des int\'er\^ets des crit\`eres \`a la Nesterenko. 

\item Les formes lin\'eaires $\ell_{i,n}-\ell_{p,n}\xi_i$ ne doivent pas \^etre trop petites dans le th\'eor\`eme \ref{Nesterenko avec diviseurs}~: on a besoin
d'une estimation exacte et pas seulement d'une majoration de $|\ell_{i,n}-\ell_{p,n}\xi_i|$. C'est l'une des diff\'erences majeures entre les crit\`eres \`a la
Nesterenko et ceux \`a la Siegel ; lorsqu'on ne fait aucune hypoth\`ese d'ind\'ependance sur
les formes lin\'eaires, une minoration de $|\ell_{i,n}-\ell_{p,n}\xi_i|$ s'av\`ere toujours n\'ecessaire (voir cependant la
proposition 1 de~\cite{MathAnnalen}),
alors qu'elle est inutile dans les crit\`eres \`a la Siegel comme le th\'eor\`eme \ref{Nesterenko a la Siegel}. 

\item Une majoration $|\ell_{p,n}|\leq Q_n^{1+o(1)}$ est suffisante dans le th\'eor\`eme \ref{Nesterenko a la Siegel} alors
qu'une \'egalit\'e est requise
dans le th\'eor\`eme \ref{Nesterenko avec diviseurs}. Cette diff\'erence est plut\^ot inhabituelle, car une majoration \`a cet endroit suffit g\'en\'eralement, m\^eme
dans les crit\`eres de type Nesterenko. 
\end{itemize}

\begin{proof}
De m\^eme que pour la d\'emonstration du th\'eor\`eme \ref{Nesterenko avec diviseurs}, on pose \break $e_i=(0,\dots,0,1,0,\dots,0,-\xi_i)$ pour $1\leq i\leq p-1$ et
$e_p=(\xi_1,\dots,\xi_{p-1},1)$. On pose $\Lambda_n=\bigoplus_{i=1}^p \delta_{i,n}\Z$ et $L_n=(\ell_{1,n},\dots,\ell_{p,n})$. On a donc $L_n\in\Lambda_n$,
$\Lambda_{n+1}\subset\Lambda_n$, $|L_n(e_i)|\leq Q_n^{-\tau_i+o(1)}$ et l'in\'egalit\'e triangulaire nous donne $\|L_n\|\leq Q_n^{1+o(1)}$. La r\'ecurrence sur les
$\ell_{i,n}$, pour tout $i\in\{1,\dots,p\}$, nous fournit la m\^eme r\'ecurrence pour la suite de formes lin\'eaires $L_n$. Et l'hypoth\`ese sur la matrice
$\Delta_{n_2}$ donne exactement l'ind\'ependance lin\'eaire de $L_{n_2},\dots,L_{n_2+p-1}$. 

On pose $\ve'=\ve/2$ et $Q'$ comme dans \eqref{Q'} de la d\'emonstration du th\'eor\`eme \ref{Nesterenko avec diviseurs}. On peut donc appliquer le th\'eor\`eme
\ref{prop a la Siegel general} et on obtient 
$$
\Lambda_{\Phi(Q')}^\perp\cap\CC_{Q',\ve'}=\{\textbf{0}\}\,.
$$
La fin de la d\'emonstration est alors similaire \`a celle du th\'eor\`eme \ref{Nesterenko avec diviseurs}.

\end{proof}


\section{Optimalit\'e}\label{Optimalite}


Dans cette section, on se place encore dans le cadre particulier qui a fait l'objet de l'\'etude du \S~\ref{Applications}. 

On v\'erifie que les conclusions des th\'eor\`emes \ref{Nesterenko avec diviseurs} et \ref{Nesterenko a la Siegel} sont optimales en construisant une
r\'eciproque \`a
ces
\'enonc\'es. Il ne s'agit cependant pas de r\'eciproques au sens propre du terme,
car les formes lin\'eaires qu'on va construire (en supposant fausses les conclusions des crit\`eres, c'est-\`a-dire les mesures d'ind\'ependance lin\'eaire restreintes)
ne satisferont pas (\emph{a priori}) \`a toutes les hypoth\`eses de ces crit\`eres. 

Dans le cas le plus simple o\`u on consid\`ere un seul nombre et pas de diviseurs, le crit\`ere de Nesterenko se r\'eduit au lemme classique suivant. 

\begin{lem}\label{lemme classique}
Soit $\xi\in\R\setminus\Q$, et $\alpha,\beta\in\R$ tels que $0<\alpha<1$ et $\beta>1$. Supposons qu'il existe des suites d'entiers $(u_n)_{n\geq 1}$ et
$(v_n)_{n\geq 1}$ tels que 

$$
\lim_{n\to+\infty}|u_n\xi-v_n|^{1/n}=\alpha\hfill\text{ et }\hfill\limsup_{n\to+\infty}|u_n|^{1/n}\leq \beta \,.
$$
Alors $\mu(\xi)\leq 1-\frac{\log \beta}{\log\alpha}$. 
\end{lem}\line

Dans ce cas particulier, Fischler et Rivoal ont quasiment d\'emontr\'e \cite{FischlerRivoal} la r\'eciproque, au sens propre du terme, de
ce r\'esultat. Ils obtiennent
m\^eme des suites $(u_n)_{n\geq 1}$ et $(v_n)_{n\geq 1}$ satisfaisant \`a des estimations asymptotiques plus pr\'ecises que celles des hypoth\`eses du lemme \ref{lemme
classique}. Le seul
d\'efaut de leur \'enonc\'e pour que ce soit v\'eritablement une r\'eciproque, est l'in\'egalit\'e stricte $\mu(\xi)<1-\frac{\log\beta}{\log\alpha}$ dans
l'hypoth\`ese, au lieu de l'in\'egalit\'e au sens large. Dans \cite{Indaga}, Fischler \'etablit \'egalement
une r\'eciproque presque compl\`ete pour l'exposant d'approximation rationnelle restreinte $\mu_{\psi}$ qui inclut des diviseurs. Ici encore, la r\'eciproque
construite n'est pas compl\`ete car l'in\'egalit\'e large du r\'esultat est prise au sens strict dans la r\'eciproque. 

Mais ces deux exemples de
r\'eciproques de crit\`eres ``\`a la Nesterenko'' ne concernent qu'une seule variable : les formes lin\'eaires construites sont des formes
lin\'eaires en 1 et $\xi$ seulement. Or, nous avons besoin de formes en plusieurs variables $1,\xi_1,\dots,\xi_{p-1}$. L'article \cite{NesterenkoConversePP} de
Fischler, Hussain, Kristensen et Levesley contient une r\'eciproque au crit\`ere de Nesterenko \`a plusieurs variables, valable presque partout au sens de la mesure
de Lebesgue. Ce r\'esultat fournit une sorte de r\'eciproque (valable
presque partout) aux th\'eor\`emes \ref{Nesterenko avec diviseurs} et \ref{Nesterenko a la Siegel} lorsque $\Lambda_n$ est le r\'eseau des formes
lin\'eaires \`a coefficients entiers. Toujours dans la direction de formes lin\'eaires en plusieurs variables,
Chantanasiri d\'emontre \'egalement (\S~3 de \cite{Chantanasiri}) une sorte de r\'eciproque
sous une hypoth\`ese tr\`es forte (une mesure d'ind\'ependance lin\'eaire). Plus de d\'etails \`a ce sujet peuvent \^etre trouv\'es dans le paragraphe 2.4 de \cite{these}.

On adopte ici la m\^eme approche en essayant d'\'etablir une sorte de r\'eciproque aux th\'eor\`emes \ref{Nesterenko avec diviseurs} et \ref{Nesterenko a la Siegel}. On montre que sous la condition
de l'existence de $p$-uplets $\mathbf a=(a_1,\dots,a_p)\in\Lambda_{\Phi(Q)}^\perp\setminus\{\textbf{0}\}$ avec $|a_i|\leq Q^{\tau_i-\ve}$ et
$|a_1\xi_1+\dots+a_{p-1}\xi_{p-1}+a_p|<Q^{-1-\ve}$ pour un $\ve>0$ donn\'e et une infinit\'e de $Q$ (ce qui est presque le contraire de la conclusion des th\'eor\`emes
\ref{Nesterenko
avec diviseurs} et \ref{Nesterenko a la Siegel}), on peut construire des entiers $\ell_{i,n}$ non tous nuls tels que $|\ell_{i,n}-\xi_i\ell_{p,n}|\leq
Q_n^{-\tau_i+o(1)}$ et $|\ell_{p,n}|\leq Q_n^{1+o(1)}$ pour $n$ assez grand (voir la proposition \ref{reciproque chap 1} ci-dessous).

On commence par montrer que les hypoth\`eses du th\'eor\`eme \ref{Nesterenko avec diviseurs} ne sont pas trop fortes et qu'il est possible de construire, pour des
r\'eels $\xi_1,\dots,\xi_{p-1}$ donn\'es, des formes lin\'eaires satisfaisant \`a la plupart des conditions du th\'eor\`eme \ref{Nesterenko avec diviseurs}. Puis, on se
pose la m\^eme question vis \`a vis des $a_i$. En effet, chacun des \'enonc\'es des parties pr\'ec\'edentes affirme qu'une forme lin\'eaire en $1,\xi_1,\dots,\xi_{p-1}$, \`a
coefficients entiers pas trop gros et suffisamment divisibles, ne peut \^etre trop petite. On d\'emontre dans la proposition \ref{construction ell et a} ci-dessous
que cette conclusion n'est pas triviale, c'est-\`a-dire que (pour certaines valeurs des param\`etres) une telle forme lin\'eaire peut \^etre tr\`es petite. 

\begin{prop}\label{construction ell et a}
Soient $\xi_1,\dots,\xi_{p-1}$ des r\'eels quelconques et $\tau_1,\dots,\tau_{p-1}$ des r\'eels $>-1$. 

On prend pour tout $i\in\{1,\dots,p\}$, et tout $n\in\N$, $\delta_{i,n}\in\N^*$ avec
$\delta_{i,n}|\delta_{i,n+1}$. On note $\Lambda_n=\delta_{1,n}\Z\oplus\dots\oplus\delta_{p,n}\Z$ le r\'eseau qu'ils
d\'efinissent. On consid\`ere aussi $(Q_n)_{n\in \N}$
une suite d'entiers strictement croissante avec $Q_{n+1}=Q_n^{1+o(1)}$. On suppose que pour tout $i\in\{1,\dots,p\}$, il
existe $\gamma_i\in\R$ tel que
$\delta_{i,n}=Q_n^{\gamma_i+o(1)}$. On note $\Phi(Q)=\max\{k\in\N, Q_k\leq Q\}$. Alors :

\begin{itemize}
\item[\textbullet] Si 
\begin{equation}\label{tau_i+gamma_i<1}
\gamma_p+\sum_{\substack{1\leq i\leq p-1 \\ \tau_i+\gamma_i\geq 0}} \tau_i+\gamma_i\leq 1\,,
\end{equation}
alors il existe, \`a partir d'un certain rang,
$(\ell_{1,n},\dots,\ell_{p,n})\in\Lambda_n\setminus\{\mathbf0\}$ tel que
$|\ell_{p,n}\xi_i-\ell_{i,n}|\leq Q_n^{-\tau_i+o(1)}$ pour $i\in\{1,\dots,p-1\}$ et $|\ell_{p,n}|\leq Q_n^{1+o(1)}$. 

\item[\textbullet] Si 
\begin{equation}\label{tau_i+gamma_i>1}
\gamma_p+\sum_{\substack{1\leq i\leq p-1 \\ \tau_i+\gamma_i\geq 0}} \tau_i+\gamma_i> 1\,,
\end{equation}
alors, pour tout $\ve>0$, tout $Q$ suffisamment grand par
rapport \`a $\ve$, il existe $(a_1,\dots,a_p)\in\Lambda_{\Phi(Q)}^\perp\setminus\{\mathbf0\}$ avec
$|a_i|\leq Q^{\tau_i-\ve}$ pour $i\in\{1,\dots,p-1\}$ et $|a_1\xi_1+\dots+a_{p-1}\xi_{p-1}+a_p|\leq Q^{-1-\ve}$.
\end{itemize}
\end{prop}

Remarquons que $\gamma_i\geq 0$ et que de plus, si $a_i\not=0$ et $Q$ est suffisamment grand, on a forc\'ement $Q_n^{-\gamma_i-\ve}\leq
\frac{1}{\delta_{i,n}}\leq |a_i|\leq Q_n^{\tau_i-\ve}$ avec $n=\Phi(Q)$, de sorte que $\gamma_i+\tau_i\geq 0$. Ainsi, $\gamma_i+\tau_i<0$ impose $a_i=0$, ce qui
motive le second point de la proposition. 

\begin{proof}
On note $J=\Big\{j\in\{1,\dots,p-1\},\tau_j+\gamma_j\geq 0\Big\}$ et $\alpha=p-|J|\in\{0,\dots,p\}$. On a donc $|J|=p-\alpha$.  On consid\`ere aussi le r\'eseau
$\Lambda'_n=\bigoplus_{j\in J\cup\{p\}}\delta_{j,n}\Z\subset\R^{J\cup\{p\}}$ qui est de d\'eterminant
$\delta_{p,n}\prod_{j\in J}\delta_{j,n}=Q_n^{\gamma_p+\sum_{j\in J}\gamma_j+o(1)}$.


Commen\c{c}ons par d\'emontrer le premier point. 

On consid\`ere $K_n=\l\{\mathbf x\in\R^{J\cup\{p\}},|x_p\xi_j-x_j|\leq Q_n^{-\tau_j+o(1)},j\in J,\text{ et } |x_p|\leq
Q_n^{1+o(1)}\r\}$. C'est un compact convexe, sym\'etrique, centr\'e en l'origine et de volume $2^{p-\alpha+1}Q_n^{1-\sum_{j\in J}\tau_j+o(1)}$. Mais l'hypoth\`ese
\eqref{tau_i+gamma_i>1} signifie
$$
1-\sum_{j\in J}\tau_j \geq \gamma_p+\sum_{j\in J}\gamma_j\,,
$$
ce qui donne :
$$
Q_n^{1-\sum_{j\in J}\tau_j+o(1)}\geq Q_n^{\gamma_p+\sum_{j\in J}\gamma_j+o(1)}\,.
$$
Le th\'eor\`eme de Minkowski assure l'existence de $\ell_{j,n}$, pour $j\in J$, et $\ell_{p,n}$ tels que $|\ell_{p,n}\xi_j-\ell_{j,n}|\leq
Q_n^{-\tau_j+o(1)}$
pour tout $j\in J$ et $|\ell_{p,n}|\leq Q_n^{1+o(1)}$.

On choisit maintenant la suite $o(1)$ comme \'etant une suite $(\ve_n)$ strictement positive telle que $\ve_n\log Q_n\not\to 0$ quand $n\to+\infty$ et on prend
un
$j\notin J$. Alors, puisque
$\tau_j+\gamma_j<0$, on a $Q_n^{\tau_j-\ve_n}\delta_{j,n}+Q_n^{-2\ve_n} \leq Q_n^{\tau_j+\gamma_j}+Q_n^{-2\ve_n}<1$ pour tout $n$ assez grand. Ainsi,
$Q_n^{-\tau_j+\ve_n}-Q_n^{-\tau_j-\ve_n}>\delta_{j,n}$ d\`es que $n$ est assez grand, et donc l'intervalle
$\l[-\ell_{p,n}\xi_j-Q_n^{-\tau_j+\ve_n},-\ell_{p,n}\xi_j-Q_n^{-\tau_j-\ve_n}\r]$ contient forc\'ement un multiple $\ell_{j,n}$ de $\delta_{j,n}$. Et
$\ell_{j,n}$ v\'erifie donc $|\ell_{j,n}-\xi_j\ell_{p,n}|\leq Q_n^{-\tau_j+\ve_n}$, ce qui finit la preuve du premier point.   

\bigskip

Montrons maintenant le deuxi\`eme point de la proposition. 

On introduit $\ve>0$ tel que 
$$
\gamma_p+\sum_{\substack{1\leq i\leq p-1 \\ \tau_i+\gamma_i\geq 0}} \tau_i+\gamma_i> 1+(p-\alpha+2)\ve\,.
$$

On consid\`ere le compact 
$$
K_Q=\l\{\mathbf a\in\R^{J\cup\{p\}},\forall j\in J, |a_j|\leq Q^{\tau_j-\ve}\text{ et }\l|\sum_{j\in J} a_j\xi_j+a_p\r|\leq Q^{-1-\ve}\r\}\,.
$$
$K_Q$ est un compact convexe,
sym\'etrique, centr\'e en l'origine, de volume $2^{p-\alpha+1}Q^{\sum_{j\in J}\tau_j-1-(p-\alpha+1)\ve}$. 

$\Lambda_n'^\perp$, le r\'eseau dual de $\Lambda_n'$, est de d\'eterminant inverse de celui de $\Lambda_n'$, c'est-\`a-dire $\prod_{j\in J\cup\{p\}}
\frac{1}{\delta_{j,n}} = Q_n^{-\gamma_p-\sum_{j\in J} \gamma_j+o(1)} = Q^{-\gamma_p-\sum_{j\in J} \gamma_j+o(1)} \leq Q^{-\gamma_p-\sum_{j\in J}
\gamma_j+\ve}$, grâce \`a la d\'efinition de $\Phi(Q)$, au comportement asymptotique de la suite $(Q_n)_{n\in\N}$ et en prenant $Q$ suffisamment grand. 

On a alors, pour $Q$ assez grand,  
$$
Q^{\sum_{j\in J}\tau_j-1-(p-\alpha+1)\ve} > Q^{-\gamma_p-\sum_{j\in J}\gamma_j +\ve} > \det(\Lambda'^\perp_{\Phi(Q)})\,.
$$

Donc $\text{Vol}(K_Q)>2^{p-\alpha+1}\det(\Lambda'^\perp_{\Phi(Q)})$ et le th\'eor\`eme de Minkowski assure l'existence d'un $(a_j)_{j\in
J\cup\{p\}}\in
K_Q\cap{\Lambda'}_{\Phi(Q)}^\perp$ non nul. On pose $a_j=0$ pour $j\not\in J\cup\{p\}$, ce qui termine la d\'emonstration. 
\end{proof}


Le deuxi\`eme point de la proposition pr\'ec\'edente permet donc de construire une forme lin\'eaire en $1,\xi_1,\dots,\xi_{p-1}$ qui contredit les conclusions des
th\'eor\`emes \ref{Nesterenko avec diviseurs} et \ref{Nesterenko a la Siegel} de la section \ref{Applications}. Ainsi, les conclusions de
ces
th\'eor\`emes impliquent que l'hypoth\`ese \eqref{tau_i+gamma_i>1} 
faite dans le second point de la proposition \ref{construction ell et a} est fausse. On peut alors appliquer le premier point et construire
des $\ell_{i,n}$ qui v\'erifient une partie des hypoth\`eses des th\'eor\`emes \ref{Nesterenko avec diviseurs} et \ref{Nesterenko a la Siegel}. Par cons\'equent,
la combinaison des deux points de la proposition \ref{construction ell et a} ci-dessus permet d'\'etablir une sorte de r\'eciproque aux th\'eor\`emes \ref{Nesterenko
avec diviseurs} et \ref{Nesterenko a la Siegel}.

\begin{prop}\label{reciproque chap 1}
Soit $\xi_1,\dots,\xi_{p-1}$ des r\'eels quelconques et $\tau_i,\dots,\tau_{p-1}$ des r\'eels $>-1$. Soit $(Q_n)_n$ une suite
d'entiers strictement croissante
telle que $Q_{n+1}=Q_n^{1+o(1)}$. Soit $p$ suites d'entiers $(\delta_{i,n})_n$ telles que pour tout $i\in\{1,\dots,p\}$, tout
$n\in\N$,
$\delta_{i,n}\big|\delta_{i,n+1}$ et
il existe $\gamma_i\in\R$ tel que $\delta_{i,n}=Q_n^{\gamma_i+o(1)}$. 

Supposons qu'il existe $\ve>0$ et une infinit\'e de $Q$ tels que pour tout $(a_1,\dots,a_p)\in\Q^p\setminus\{\mathbf0\}$ tel que
$(\delta_{1,\Phi(Q)}a_1,\dots,\delta_{p,\Phi(Q)}a_p)\in\Z^p$ avec $|a_i|\leq Q^{\tau_i-\ve}$ pour $1\leq i\leq p-1$, on
ait
$|a_1\xi_1+\dots+a_{p-1}\xi_{p-1}+a_p|>Q^{-1-\ve}$. 

Alors pour tout $n$ assez grand, il existe $\ell_{i,n}\in\delta_{i,n}\Z$ pour $i\in\{1,\dots,p\}$, non tous nuls, 
tels que pour $1\leq i\leq p-1$, $|\ell_{i,n}-\xi_i\ell_{p,n}|\leq Q_n^{-\tau_i+o(1)}$ et $|\ell_{p,n}|\leq Q_n^{1+o(1)}$. 
\end{prop}

La double absence de relation de r\'ecurrence et d'ind\'ependance lin\'eaire entre des rangs successifs des formes lin\'eaires construites fait que ce n'est pas exacteme,t la r\'eciproque du crit\`ere \`a la Siegel du th\'eor\`eme \ref{Nesterenko a la Siegel}. Ce n'est pas une r\'eciproque exacte du 
crit\`ere \`a la Nesterenko du th\'eor\`eme \ref{Nesterenko avec diviseurs}, les
estimations sur les tailles des formes lin\'eaires et de leurs coefficients sont exactes, alors que nous ne sommes en mesure ici d'obtenir que des majorations.

\bibliographystyle{alpha-fr}
\bibliography{Biblio_Theorie}

\begin{thebibliography}{{Mar}06}
\expandafter\ifx\csname fonteauteurs\endcsname\relax
\def\fonteauteurs{\scshape}\fi

\bibitem[Bed98]{Bedulev}
Egor~V. \bgroup\fonteauteurs\bgroup Bedulev\egroup\egroup{} :
\newblock On the linear independence of numbers over number fields.
\newblock {\em Mat. Zametki}, 64(4)\string:\penalty500\relax 506--517, 1998.

\bibitem[BR01]{IrrationaliteInfValeursZeta}
Keith \bgroup\fonteauteurs\bgroup Ball\egroup\egroup{} et Tanguy
  \bgroup\fonteauteurs\bgroup Rivoal\egroup\egroup{} :
\newblock Irrationalit\'e d'une infinit\'e de valeurs de la fonction z\^eta aux
  entiers impairs.
\newblock {\em Invent. Math.}, 146(1)\string:\penalty500\relax 193--207, 2001.

\bibitem[Cas97]{Cassels}
John W.~S. \bgroup\fonteauteurs\bgroup Cassels\egroup\egroup{} :
\newblock {\em An introduction to the geometry of numbers}.
\newblock Classics in {M}athematics. Springer-{V}erlag, Berlin, 1997.
\newblock Corrected reprint of the 1971 edition.

\bibitem[Cha12]{Chantanasiri}
Amarisa \bgroup\fonteauteurs\bgroup Chantanasiri\egroup\egroup{} :
\newblock G\'en\'eralisation des crit\`eres pour l'ind\'ependance lin\'eaire de
  {N}esterenko, {A}moroso, {C}olmez, {F}ischler et {Z}udilin.
\newblock {\em Ann. Math. Blaise Pascal}, 19(1)\string:\penalty500\relax
  75--105, 2012.

\bibitem[Dau14]{these}
Simon \bgroup\fonteauteurs\bgroup Dauguet\egroup\egroup{} :
\newblock {\em Généralisations du critère d'indépendance linéaire de
  Nesterenko}.
\newblock Th\`ese de doctorat, Laboratoire de Mathématiques d'Orsay,
  Paris-Sud, 2014.
\newblock \url{http://www.math.u-psud.fr/~dauguet/These.html}.

\bibitem[DZ14]{DauZud}
Simon \bgroup\fonteauteurs\bgroup Dauguet\egroup\egroup{} et Wadim
  \bgroup\fonteauteurs\bgroup Zudilin\egroup\egroup{} :
\newblock On simultaneous diophantine approximations to $\zeta(2)$ and
  $\zeta(3)$.
\newblock \href{http://arxiv.org/abs/1401.5322}{arXiv:1401.5322}. Soumis, 2014.

\bibitem[FHKL13]{NesterenkoConversePP}
Stéphane \bgroup\fonteauteurs\bgroup Fischler\egroup\egroup{}, Mumtaz
  \bgroup\fonteauteurs\bgroup Hussain\egroup\egroup{}, Simon
  \bgroup\fonteauteurs\bgroup Kristensen\egroup\egroup{} et Jason
  \bgroup\fonteauteurs\bgroup Levesley\egroup\egroup{} :
\newblock A converse to linear independence criteria, valid almost everywhere.
\newblock \href{http://arxiv.org/abs/1302.1952}{arXiv:1302.1952}, Soumis, 2013.

\bibitem[Fis09]{Indaga}
St\'ephane \bgroup\fonteauteurs\bgroup Fischler\egroup\egroup{} :
\newblock Restricted rational approximation and {A}p\'ery-type constructions.
\newblock {\em Indagationes Mathematicae}, 20 (2)\string:\penalty500\relax
  201--215, Juin 2009.

\bibitem[Fis12]{NesterenkoOscillate}
St{\'e}phane \bgroup\fonteauteurs\bgroup Fischler\egroup\egroup{} :
\newblock Nesterenko's criterion when the small linear forms oscillate.
\newblock {\em Arch. Math. (Basel)}, 98(2)\string:\penalty500\relax 143--151,
  2012.

\bibitem[Fis13]{Nesterenko4Vectors}
Stéphane \bgroup\fonteauteurs\bgroup Fischler\egroup\egroup{} :
\newblock Nesterenko's linear independence criterion for vectors.
\newblock \href{http://arxiv.org/abs/1202.2279v2}{arXiv:1202.2279v2}. Soumis,
  Octobre 2013.

\bibitem[FN98]{NumberTheory}
Naum~I. \bgroup\fonteauteurs\bgroup Fel'dman\egroup\egroup{} et Yuri~V.
  \bgroup\fonteauteurs\bgroup Nesterenko\egroup\egroup{} :
\newblock Transcendental numbers.
\newblock \emph{In} {\em Number theory, {IV}}, volume~44 de {\em Encyclopaedia
  Math. Sci.}, pages 1--345. Springer, Berlin, 1998.

\bibitem[FR10]{FischlerRivoal}
St{\'e}phane \bgroup\fonteauteurs\bgroup Fischler\egroup\egroup{} et Tanguy
  \bgroup\fonteauteurs\bgroup Rivoal\egroup\egroup{} :
\newblock Irrationality exponent and rational approximations with prescribed
  growth.
\newblock {\em Proc. Amer. Math. Soc.}, 138(3)\string:\penalty500\relax
  799--808, 2010.

\bibitem[FZ10]{MathAnnalen}
Stéphane \bgroup\fonteauteurs\bgroup Fischler\egroup\egroup{} et Wadim
  \bgroup\fonteauteurs\bgroup Zudilin\egroup\egroup{} :
\newblock {A} refinement of {N}esterenko's linear independence criterion with
  applications to zeta values.
\newblock {\em Math. Annalen}, 347\string:\penalty500\relax 739--763, 2010.

\bibitem[{Mar}06]{Marcovecchio}
Raffaele \bgroup\fonteauteurs\bgroup {Marcovecchio}\egroup\egroup{} :
\newblock {Linear independence of linear forms in polylogarithms.}
\newblock {\em {Ann. Sc. Norm. Super. Pisa, Cl. Sci. (5)}},
  5(1)\string:\penalty500\relax 1--11, 2006.

\bibitem[Nes85]{Nesterenko}
Yuri~V. \bgroup\fonteauteurs\bgroup Nesterenko\egroup\egroup{} :
\newblock Linear independence of numbers.
\newblock {\em Vestnik Moskov. Univ. Ser. I Mat. Mekh.},
  (1)\string:\penalty500\relax 46--49, 108, 1985.

\bibitem[Nes12]{Nesterenkopadic}
Yuri~V. \bgroup\fonteauteurs\bgroup Nesterenko\egroup\egroup{} :
\newblock On a criterion of linear independence of {$p$}-adic numbers.
\newblock {\em Manuscripta Math.}, 139(3-4)\string:\penalty500\relax 405--414,
  2012.

\bibitem[Riv00]{FctZetaInfValeursIrrationnelles}
Tanguy \bgroup\fonteauteurs\bgroup Rivoal\egroup\egroup{} :
\newblock La fonction zêta de {R}iemann prend une infinité de valeurs
  irrationnelles aux entiers impairs.
\newblock {\em C.R.A.S. Paris Série I Math.}, 336.4\string:\penalty500\relax
  267--270, 2000.

\bibitem[Riv02]{Rivoal21}
Tanguy \bgroup\fonteauteurs\bgroup Rivoal\egroup\egroup{} :
\newblock Irrationalit\'e d'au moins un des neuf nombres
  {$\zeta(5),\zeta(7),\dots,\zeta(21)$}.
\newblock {\em Acta Arith.}, 103(2)\string:\penalty500\relax 157--167, 2002.

\bibitem[Sie29]{Siegel}
Carl~L. \bgroup\fonteauteurs\bgroup Siegel\egroup\egroup{} :
\newblock {Ü}ber einige {A}nwendungen {D}iophantischer {A}pproximationen.
\newblock {\em Abh. Press. Akad. Wiss. Phys.-Math. Kl.},
  1\string:\penalty500\relax 1--70, 1929.
\newblock JFM 56.0180.01.

\bibitem[T{\"o}p94]{Topfer1}
Thomas \bgroup\fonteauteurs\bgroup T{\"o}pfer\egroup\egroup{} :
\newblock An axiomatization of {N}esterenko's method and applications on
  {M}ahler functions.
\newblock {\em J. Number Theory}, 49(1)\string:\penalty500\relax 1--26, 1994.

\bibitem[T{\"o}p95]{Topfer2}
Thomas \bgroup\fonteauteurs\bgroup T{\"o}pfer\egroup\egroup{} :
\newblock An axiomatization of {N}esterenko's method and applications on
  {M}ahler functions. {II}.
\newblock {\em Compositio Math.}, 95(3)\string:\penalty500\relax 323--342,
  1995.

\bibitem[Zud01]{Irrationalite4}
Wadim \bgroup\fonteauteurs\bgroup Zudilin\egroup\egroup{} :
\newblock One of the numbers $\zeta(5)$, $\zeta(7)$, $\zeta(9)$, $\zeta(11)$ is
  irrational.
\newblock {\em Russian Math. Surveys}, 56.4\string:\penalty500\relax 774--776,
  2001.

\bibitem[Zud02]{Irr145}
Wadim \bgroup\fonteauteurs\bgroup Zudilin\egroup\egroup{} :
\newblock Irrationality of values of the {R}iemann zeta function.
\newblock {\em Izv. Math.}, 3\string:\penalty500\relax 489--542, 2002.

\end{thebibliography}

\end{document}